\documentclass[journal]{new-aiaa}

\usepackage{amsthm}
\usepackage{amsmath}
\usepackage{amsfonts}
\usepackage{hyperref}
\usepackage{graphicx}
\usepackage{subcaption}
\usepackage{blindtext}
\usepackage{epsfig} 

\newtheorem{theorem}{Theorem}[section]

\newtheorem{prop}[theorem]{\textbf{Proposition}}

\usepackage{algorithmic}
\usepackage[ruled,linesnumbered]{algorithm2e}
\usepackage{bm}

\title{Learning Linear Time-Varying Koopman Models for Constrained Nonlinear Optimal Control}

\author{ Zhi Xu \footnote{Graduate Research Assistant, School of Aeronautics and Astronautics.  {\tt\small xu1710@purdue.edu}}}
\author{Ran Dai \footnote{Professor, School of Aeronautics and Astronautics. AIAA Associate Fellow. Corresponding Author, {\tt\small randai@purdue.edu}}
}
\affil{Purdue University, West Lafayette, IN, 47907}

\begin{document}
\footnotetext{Presented as Paper 2026-1168 at the AIAA Scitech 2026 Forum, Orlando, FL, January 12–16, 2026}

\maketitle

\begin{abstract}
{
This paper presents a learning-based Koopman framework for solving constrained nonlinear optimal control problems. Despite significant advances in Koopman-based modeling and control, challenges remain in achieving accurate representations of controlled nonlinear dynamics and constraints while maintaining computational efficiency for optimal control. Motivated by these challenges, a linear time-varying (LTV) Koopman representation is developed to provide a more expressive representation of controlled nonlinear dynamics while preserving a linear prediction structure suitable for efficient convex optimization. Modeling errors and their propagation over extended prediction horizons are analyzed, motivating a neural-network-based Koopman framework that jointly learns lifted dynamics and constraint representations with stability-oriented regularization to improve long-horizon prediction performance. Based on the learned Koopman model, the original constrained nonlinear optimal control problem is reformulated as a convex quadratic program with conservative constraint tightening to account for modeling uncertainty during constraint enforcement. Numerical results on a low-thrust orbit transfer problem and a Mars entry guidance problem demonstrate the effectiveness of the proposed framework in terms of prediction accuracy, computational efficiency, and constrained optimal control performance.
}
\end{abstract}

\section{Introduction}

Nonlinear dynamical systems subject to state and input constraints arise ubiquitously in engineering and scientific applications. In many practical settings, the system dynamics exhibit strong nonlinearities, while safety, performance, or physical limitations impose nonlinear constraints that must be satisfied along the entire trajectory. The simultaneous presence of nonlinear dynamics and nonlinear constraints poses significant challenges for system analysis, prediction, and control. Although nonlinear models can capture the underlying physics with high fidelity, they typically lead to nonconvex optimal control problems that are difficult to solve efficiently, particularly when constraints must be enforced explicitly. As a result, there remains a fundamental need for modeling and control frameworks that can handle nonlinear dynamics and constraints in a computationally tractable manner.

A wide range of approaches has been developed to address nonlinear constrained optimal control problems (OCPs), broadly categorized into direct and indirect methods. Indirect methods derive necessary conditions for optimality through Pontryagin’s minimum principle, resulting in two-point boundary value problems involving both states and costates \cite{Bryson1979, Betts1998SurveyON, Zheng2017I}. While such methods can yield highly accurate solutions, their practical application is often hindered by sensitivity to the initial guesses of costates, whose lack of direct physical interpretation makes reliable initialization difficult. Moreover, the incorporation of complex state and path constraints within indirect formulations remains challenging \cite{arutyunov2020survey}. Direct methods, on the other hand, transcribe the nonlinear OCP into a finite-dimensional nonlinear program (NLP) using discretization techniques, e.g., collocation \cite{hargraves1987direct, gong2009pseudospectral, liu2025survey}. Direct methods offer greater flexibility in handling constraints and have therefore become the dominant paradigm in practice. However, the resulting optimization problems are typically large-scale and nonconvex, with computational cost and convergence behavior that depend critically on discretization resolution and solver initialization \cite{schimke2023rapid, Liu2022M}. In addition, the nonlinear transcription may lead to convergence to locally optimal solutions. The substantial computational burden associated with solving such NLPs further limits their applications in resource-constrained or real-time settings \cite{liu2017survey}.

Motivated by the limitations of solving nonlinear OCPs, alternative modeling frameworks have been explored to balance expressiveness and computational tractability. Among these, Koopman operator theory has attracted increasing attention as a powerful approach for representing nonlinear dynamical systems through linear evolution of observables in the infinite-dimensional function space \cite{Brunton2022MKT}. By transferring nonlinear dynamics into a linear form, Koopman-based representations offer a promising foundation for analysis and control using tools from linear systems theory. However, practical application requires constructing a finite-dimensional approximation of the Koopman operator. For example, the dynamic mode decomposition (DMD) algorithm approximates the unitary Koopman operator using linear observables \cite{Schmid2008D}, but its expressiveness is limited by the linear structure. The extended dynamic mode decomposition (EDMD) method addresses this limitation by augmenting the observable set with nonlinear functions, thereby creating a higher-dimensional lifted representation with improved nonlinear modeling capability \cite{Williams2014A}. More recently, neural-network–based techniques have been introduced to learn the observable dictionary from data, leveraging their expressive power to produce more accurate approximations \cite{Takeishi2017L, Li2017E, Lusch2018}. {
Despite these advances, the application of Koopman-based methods to constrained nonlinear optimal control remains challenging. In particular, accurately representing controlled nonlinear dynamics while preserving computational tractability, and systematically handling nonlinear state and path constraints within the lifted space, remain important open problems.

One of the primary challenges arises from extending the Koopman framework to controlled systems. {\color{black} Although the Koopman operator theory is well established for nonlinear dynamical systems, its extension to controlled systems remains challenging \cite{Otto2021K}. In general, the presence of control inputs breaks the semigroup property that underlies the Koopman operator for autonomous systems, impeding both theoretical analysis and practical approximation. Existing Koopman-based control approaches often adopt a linear time-invariant (LTI) predictor structure with either fixed lifting functions~\cite{Korda2018L} or neural-network-based lifting functions~\cite{han2020deep}. Such formulations implicitly impose a linear dependence on the control variables in the lifted dynamics, making it difficult to accurately represent general nonlinear control effects.} In addition, the use of a single invariant operator may become overly restrictive for systems exhibiting significant temporal variations or large excursions across the state space. Consequently, maintaining prediction accuracy over a broad operating envelope often requires either a very rich lifting space or a more flexible operator representation.

{\color{black}To better capture the control-dependent structure of nonlinear systems, bilinear Koopman representations have been developed for control-affine systems, where the Koopman generator admits an exact bilinear decomposition with respect to the control inputs \cite{Zhao2024, SafEDMD2026}.} This formulation preserves the control-affine structure of the original nonlinear system and provides a more faithful representation of state–control interactions. However, the resulting lifted dynamics remain bilinear in the state and control variables, and the associated optimal control problems generally remain nonconvex. Consequently, many of the computational advantages of linear prediction and convex optimization cannot be directly exploited in practical control synthesis when applying the bilinear representation.

A further challenge lies in the treatment of nonlinear state and path constraints within Koopman-based optimal control formulations. {\color{black}Conventional strategies often rely on manually designed observable dictionaries or heuristic constraint lifting strategies to obtain approximately linear constraint representations in the lifted space \cite{Korda2018L,Mauroy2020I,Otto2021K}.} For complex nonlinear constraints, however, it is often unclear whether a suitable set of observables exists that can accurately represent the constraint while remaining compatible with the lifted dynamics. Moreover, approximation errors introduced during constraint lifting may lead to a mismatch between constraint satisfaction in the lifted space and constraint satisfaction for the original nonlinear system, limiting the reliability of Koopman-based constrained optimal control.

To address the limitations of existing Koopman-based approaches for constrained optimal control, this work develops a learning-based linear time-varying Koopman framework that seeks to balance modeling fidelity, computational tractability, and constraint reliability. The central idea is to construct a linear time-varying (LTV) Koopman representation that retains a linear prediction structure while relaxing the restrictive assumption of a single invariant Koopman operator. By allowing the lifted dynamics to evolve over time, the proposed representation provides greater flexibility for capturing nonlinear system behaviors across broad operating conditions and long prediction horizons. Building upon this representation, we develop a data-driven framework that jointly learns lifted dynamics and constraint representations. To improve reliability in long-horizon prediction and constrained control, approximation errors and their propagation are analyzed and incorporated into both the learning and control design. The resulting Koopman model enables the original constrained nonlinear optimal control problem to be reformulated as a computationally tractable convex optimization problem with explicit consideration of modeling uncertainty.

The main contributions of this work can be summarized as follows:
\begin{itemize}
    \item We develop an LTV Koopman representation for controlled nonlinear systems by constructing an LTV approximation of a discrete-time bilinear Koopman model. The resulting representation preserves a linear prediction structure while providing greater flexibility than conventional LTI Koopman models.
    \item We propose a neural-network-based Koopman learning framework that jointly learns lifted dynamics and constraint representations from data. Motivated by the approximation error analysis, a stability-oriented regularization strategy is incorporated to improve long-horizon prediction performance.
    \item We formulate a Koopman-based convex optimal control framework for constrained nonlinear systems. An error-aware constraint-tightening strategy is introduced to account for approximation errors when enforcing lifted-space constraints.
\end{itemize}
}

The remainder of the paper is organized as follows. Section \ref{sec:LTV} presents the proposed LTV Koopman representation and its associated error analysis. Section \ref{sec:NNKF} introduces the neural-network-based framework for jointly learning lifted dynamics and constraint representations. Section \ref{sec:koopman_oc} develops the Koopman-based convex optimal control formulation with error-aware constraint tightening. {Section \ref{sec:simulation} presents numerical results for a low-thrust orbit transfer problem and a six-degree-of-freedom (6-DOF) Mars entry guidance problem. Finally, Section \ref{sec:conclusion} concludes the paper.}

{
\section{LTV Koopman Representation and Error Analysis}
\label{sec:LTV}
This section establishes the theoretical foundation of the proposed framework. Starting from a control-affine nonlinear system, we first derive a bilinear Koopman representation and its corresponding discrete-time approximation. Building upon this model, a linear time-varying (LTV) Koopman formulation is constructed to obtain a representation that is compatible with linear prediction and optimization frameworks. Finally, approximation errors introduced by discretization and bilinear-to-linear approximation are analyzed, and their propagation over long prediction horizons is characterized. These results provide the basis for the learning and optimization frameworks developed in the subsequent sections.

\subsection{Koopman Representation of Control-Affine Systems}
Consider a continuous-time control-affine nonlinear system
\begin{equation}\label{eq:affine}
\dot{\mathbf{x}}(t)
=
\mathbf{f}(\mathbf{x}(t))
+
\sum_{i=1}^{m}
\mathbf{g}_i(\mathbf{x}(t))
\,\mathbf{u}_i(t),
\quad
\mathbf{x}(t)\in\mathbb{R}^{n},
\quad
\mathbf{u}(t)\in\mathbb{R}^{m},
\end{equation}
where state $\mathbf{x}\in\mathbb{R}^n$, control $\mathbf{u}\in\mathbb{R}^m$, dynamics $\mathbf{f}, \mathbf{g}_i$ are locally Lipschitz.
Instead of studying the evolution of the state directly, Koopman operator theory characterizes the dynamics through observable functions. Let
$\psi:\mathbb{R}^{n}\rightarrow\mathbb{R}$
be a continuously differentiable observable.

For autonomous systems, the Koopman operator acts on observables by composition with the flow map and forms a linear semigroup. For control-affine systems, however, the dynamics are non-autonomous and the classical Koopman semigroup property does not directly apply. Nevertheless, the affine structure of the control input allows the infinitesimal evolution of observables to be decomposed into drift and control-dependent components.

Applying the chain rule along system trajectories yields
\[
\frac{d}{dt}\psi(\mathbf{x}(t))
=
\nabla\psi(\mathbf{x}(t))
\left(
\mathbf{f}(\mathbf{x}(t))
+
\sum_{i=1}^{m}
\mathbf{g}_{i}(\mathbf{x}(t))
\,\mathbf{u}_{i}(t)
\right).
\]
Define the Koopman generators associated with the drift and control vector fields as
\[
[\mathcal{K}_{0}\psi](\mathbf{x})
:=
\nabla\psi(\mathbf{x})\mathbf{f}(\mathbf{x}),
\qquad
[\mathcal{K}_{i}\psi](\mathbf{x})
:=
\nabla\psi(\mathbf{x})\mathbf{g}_{i}(\mathbf{x}),
\quad
i=1,\ldots,m.
\]
Then the observable dynamics can be expressed as
\[
\frac{d}{dt}\psi(\mathbf{x}(t))
=
[\mathcal{K}_{0}\psi](\mathbf{x}(t))
+
\sum_{i=1}^{m}
[\mathcal{K}_{i}\psi](\mathbf{x}(t))
\,\mathbf{u}_{i}(t).
\]

Assume there exists a finite-dimensional lifting function
$
\Psi(\mathbf{x}) = [\psi_1(\mathbf{x})\ \psi_2(\mathbf{x})\ \cdots\ \psi_N(\mathbf{x})]^\top
$
whose span is invariant under $\mathcal{K}_0,\mathcal{K}_1,\dots,\mathcal{K}_m$. Then there exist matrices $K_0,K_1,\dots,K_m\in\mathbb{R}^{N\times N}$ such that
\[
\frac{d}{dt}\Psi(\mathbf{x}) = K_0\Psi(\mathbf{x}) + \sum_{i=1}^m u_i K_i \Psi(\mathbf{x}).
\]
Defining the lifted state $\mathbf{z}=\Psi(\mathbf{x})$, the lifted dynamics can then be expressed as
\begin{equation} \label{eq:bilinear_dyn}
\dot{\mathbf{z}} = K_0\mathbf{z} + \sum_{i=1}^m u_i K_i \mathbf{z},
\end{equation}
which is a continuous-time bilinear system in the lifted Koopman coordinates.

{\color{black}To obtain a tractable discrete-time representation for controller design, we integrate the lifted dynamics~\eqref{eq:bilinear_dyn} over a sampling interval $\Delta t > 0$, yielding
\[
\mathbf{z}_{k+1}
=
\mathbf{z}(t_k+\Delta t)
=
\mathbf{z}(t_k)
+
\int_{t_k}^{t_k+\Delta t}
\dot{\mathbf{z}}(t)\,dt.
\]
We assume that the control input is piecewise constant within each sampling interval, such that
\[
\mathbf{u}(t)
=
\mathbf{u}_{k},
\qquad
t\in[k\Delta t,(k+1)\Delta t).
\]
Substituting the continuous-time lifted dynamics into the integral equation gives
\begin{equation} \label{eq:affine_discrete}
\mathbf{z}_{k+1}
=
\mathbf{z}_{k}
+
\int_{0}^{\Delta t}
\left(
K_{0}\mathbf{z}(t_k+\tau)
+
\sum_{i=1}^{m}
K_{i}\mathbf{z}(t_k+\tau)
\,\mathbf{u}_{k,i}
\right)d\tau .
\end{equation}

By employing a first-order Taylor expansion to approximate the integrand at $\tau = 0$, the discrete-time update can be expressed as
\[
\mathbf{z}_{k+1}
=
\mathbf{z}_{k}
+
\Delta t
\left(
K_{0}\mathbf{z}_{k}
+
\sum_{i=1}^{m}
K_{i}\mathbf{z}_{k}
\,\mathbf{u}_{k,i}
\right)
+
O(\Delta t^{2}).
\]
Finally, by defining the discrete-time system matrices as
\[
A
:=
I+\Delta t\,K_{0},
\qquad
B_i
:=
\Delta t\,K_{i},
\]
we arrive at the first-order accurate discrete-time bilinear Koopman model of the continuous-time control-affine system:
\begin{equation} \label{eq:bkm}
\mathbf{z}_{k+1}
=
A\mathbf{z}_{k}
+
\sum_{i=1}^{m}
B_i\mathbf{z}_{k}\mathbf{u}_{k,i}
+
O(\Delta t^{2}).
\end{equation}}
The resulting discrete-time bilinear Koopman model~\eqref{eq:bkm} serves as an intermediate representation between the original nonlinear dynamics and the linear time-varying model. Compared with conventional LTI Koopman formulations, the bilinear representation preserves the control-affine structure of the original system and provides a more expressive description of state-control interactions. The discretization remainder term ($O(\Delta t^2)$) explicitly characterizes the approximation introduced from discretization and constitutes one source of prediction error. This representation forms the foundation for the subsequent construction of an LTV Koopman model and the associated error analysis.
}

{
\subsection{LTV Koopman Approximation}
While the bilinear Koopman model provides a faithful representation of control-affine dynamics, its bilinear dependence on the control inputs prevents the direct application of linear prediction and optimization techniques. In particular, the resulting optimal control problem adopting bilinear dynamics generally remains nonconvex, limiting the computational advantages typically associated with Koopman-based control. To obtain a representation that is compatible with standard linear control and optimization frameworks, we construct a linear time-varying (LTV) approximation of the bilinear Koopman model. The approximation error introduced by this transformation is then quantified to assess its impact on prediction accuracy.

\begin{theorem}{\color{black}[Linear Time-Varying Koopman Model]
Consider the linear time-varying Koopman model
\begin{equation} \label{eq:ltv}
\hat{\mathbf z}_{k+1}
=
A_k \mathbf z_k + B_k \mathbf u_k + \mathbf d_k .
\end{equation}}
Define the one-step approximation error as
\[
\varepsilon_k
:=
\mathbf z_{k+1}-\hat{\mathbf z}_{k+1},
\]
where $\mathbf z_{k+1}$ denotes the exact lifted state generated by the control-affine update~\eqref{eq:affine_discrete}.
Then
{\color{black}\[
\|\varepsilon_k\|
\le
C_1\Delta t^2+C_2\Delta t,
\]
for some constants $C_1,C_2>0$.}
\end{theorem}
}
\begin{proof}
We start from the discrete-time bilinear Koopman model
\begin{equation} \label{eq:bkm_dyn}
 \mathbf{z}_{k+1} = A  \mathbf{z}_k + \sum_{i=1}^m B_i  \mathbf{z}_k \mathbf{u}_{k,i} +\mathbf{\varepsilon}_d,
\end{equation}
where $A = I + \Delta t\,{K}_0$ and $B_i = \Delta t\,{K}_i$ are the discrete-time Koopman operators, the discretization error $\varepsilon_d$ satisfies
$
\|\varepsilon_d\| \le C_1 \Delta t^2
$ for some constant $C_1>0$. 

By defining a linear time-varying operator
\[
K_k := A + \sum_{i=1}^m B_i  \mathbf{u}_{k,i},
\]
the dynamics~\eqref{eq:bkm_dyn} can be rewritten as $\mathbf{z}_{k+1} = K_k  \mathbf{z}_k + \varepsilon_d$.

To obtain a standard linear control input form, we decompose the bilinear term around a nominal point ($\bar{\mathbf{z}}_k, \bar{\mathbf{u}}_k$):
\begin{align} \nonumber
\mathbf{z}_{k+1} &= A \mathbf{z}_k + \sum_{i=1}^m B_i ( \bar{\mathbf{z}}_k \mathbf{u}_{k,i} + \mathbf{z}_k \bar{\mathbf{u}}_{k,i} - \bar{\mathbf{z}}_k \bar{\mathbf{u}}_{k,i} + \delta \mathbf{z}_k \delta \mathbf{u}_{k,i} ) + \varepsilon_d, \\  \label{eq:ltv_proof}
&= A_k \mathbf{z}_k + B_k \mathbf{u}_k + \mathbf{d}_k + \varepsilon_l + \varepsilon_d
\end{align}
where the linear time-varying operators are defined as
\[
A_k := A + \sum_{i=1}^m B_i \bar{\mathbf{u}}_{k,i},\quad 
B_k := \begin{bmatrix} B_1 \bar { \mathbf{z}}_k & \cdots & B_m \bar { \mathbf{z}}_k \end{bmatrix},\quad
\mathbf{d}_k := - \sum_{i=1}^m B_i\bar{\mathbf{z}}_k \bar{\mathbf{u}}_{k,i},
\]
and the residual term
\[
\varepsilon_l := \sum_{i=1}^m B_i\, \delta \mathbf{z}_k\, \delta \mathbf{u}_{k,i}
\]
represents the second-order error introduced by bilinear deviations, where 
$\delta \mathbf{z}_k = \mathbf{z}_k - \bar{\mathbf{z}}_k,\
\delta \mathbf{u}_{k,i} = \mathbf{u}_{k,i} - \bar{\mathbf{u}}_{k,i}$.

Assume that the deviations $\delta \mathbf{z}_k$ and $\delta \mathbf{u}_k$ are uniformly bounded, i.e., $\|\delta \mathbf{z}_k\| \le c_z$ and $\|\delta \mathbf{u}_k\| \le c_u$.
Using a compatible induced matrix norm and recalling that $B_i = \Delta t\,{K}_i$, we obtain
\[
\|\varepsilon_l\|
\le \sum_{i=1}^m \|B_i\| \, \|\delta \mathbf{z}_k\| \, \|\delta \mathbf{u}_{k,i}\|
\le \Delta t \, c_z c_u \sum_{i=1}^m \|{K}_i\|
= C_2 \Delta t,
\]
where the constant $C_2 = c_z c_u \sum_{i=1}^m \|{K}_i\|$ is independent of $\Delta t$.

The total one-step approximation error is given by
$\varepsilon_k = \varepsilon_d + \varepsilon_l.$
Applying the triangle inequality yields
\[
\|\varepsilon_k\|
\le
\|\varepsilon_d\|+\|\varepsilon_l\|
\le
C_1\Delta t^2+C_2\Delta t .
\]
Substituting $\varepsilon_k$ into equation~\eqref{eq:ltv_proof} and defining $\hat{\mathbf z}_{k+1}
=
A_k \mathbf z_k + B_k \mathbf u_k + \mathbf d_k$, we obtain the linear time-varying model~\eqref{eq:ltv} with bounded one-step approximation error $\varepsilon_k$.
\end{proof}

For dynamics with nonlinear control dependence, a first-order Taylor expansion with respect to the input about $\bar{\mathbf{u}}_k$ yields a control-affine approximation with a quadratic remainder. Under the same lifting assumption for this approximation, and assuming bounded second derivatives with respect to the input and bounded lifting derivatives over the operating domain, this remainder introduces an additional one-step error of order $O(\Delta t\|\delta\mathbf{u}_k\|^2)$, beyond the discretization and bilinear approximation errors considered above.

\subsection{Multi-step Error Propagation Analysis}\label{sec:error propagation}
{
The error derived in the previous subsection characterizes the approximation accuracy of the proposed LTV Koopman model within the one-step horizon. However, optimal control applications typically require repeated propagation of the approximation over extended prediction horizons. Therefore, it is important to understand how one-step modeling errors propagate through the time-varying Koopman dynamics and how the resulting prediction error depends on the properties of operators.

\begin{prop}[Multi-step Prediction Error]
Let $T=M\Delta t$ be a prediction horizon consisting of $M$ discrete-time steps.
Consider the $M$-step prediction $\hat{\mathbf{z}}_M$ {\color{black}generated by the linear time-varying Koopman model
\[
\hat{\mathbf z}_{k+1}
=
A_k \hat{\mathbf z}_k
+
B_k \mathbf u_k
+
\mathbf d_k,
\qquad
\hat{\mathbf z}_0=\mathbf z_0 .
\]}
Define the multi-step prediction error at horizon $T$ as
\[
\mathbf e_M
:=
\mathbf z_M-\hat{\mathbf z}_M,
\]
where $\mathbf z_M$ denotes the exact lifted state generated by the control-affine update~\eqref{eq:affine_discrete}.
Then
{\color{black}\begin{equation}
\|\mathbf{e}_M\| \le \frac{(1 + \alpha \Delta t)^M - 1}{\alpha} ( C_1 \Delta t + C_2),
\end{equation}
where $\alpha > 0$ is a constant related to the induced norm of $A_k$.}
\end{prop}
}

\begin{proof}
The prediction error at step $k+1$ can be written as
\begin{equation}
\mathbf{e}_{k+1} := \mathbf{z}_{k+1} - \hat{\mathbf{z}}_{k+1} = A_k \mathbf{e}_k + \varepsilon_k,
\end{equation}
where $\mathbf{e}_k = \mathbf{z}_k - \hat{\mathbf{z}}_k$.

Expanding this recurrence from the initial error $\mathbf{e}_0 = \mathbf{0}$, and applying the triangle inequality and induced matrix norm, we obtain
\begin{equation}
\|\mathbf{e}_M\| = \|\sum_{k=0}^{M-1} ( \prod_{j=k+1}^{M-1} A_j ) \varepsilon_k\| \le \sum_{k=0}^{M-1} ( \prod_{j=k+1}^{M-1} \|A_j\| ) \|\varepsilon_k\|.
\end{equation}
Recall that
\begin{equation}
A_j = I + \Delta t ({K}_0 + \sum_{i=1}^m {K}_i \bar{\mathbf{u}}_{j,i})
\end{equation}
Assume that the Koopman generators are bounded and that the nominal input is uniformly bounded, then there exists a constant $\alpha > 0$ such that
\begin{equation}
\|A_j\| \le 1 + \alpha \Delta t.
\end{equation}
Substituting the single-step error bound and simplifying it using geometric series, we arrive at the multi-step error bound
\begin{equation}
\|\mathbf{e}_M\| \le \sum_{k=0}^{M-1} (1 + \alpha \Delta t)^{M-1-k} ( C_1 \Delta t^2 + C_2 \Delta t ) = \frac{(1 + \alpha \Delta t)^M - 1}{\alpha} ( C_1 \Delta t + C_2 ).
\end{equation}
\end{proof}

To analyze the behavior of the error bound as the sampling frequency increases, we next consider the limit $\Delta t \to 0$ while keeping the total time horizon $T = M \Delta t$ fixed. Using the fundamental limit, we have:
\begin{equation}
\lim_{\Delta t \to 0} (1 + \alpha \Delta t)^M = \lim_{\Delta t \to 0} (1 + \alpha \Delta t)^{T/\Delta t} = e^{\alpha T}.
\end{equation}
Consequently, the continuous-time error evolution satisfies:
\begin{equation}
\|\mathbf{e}(T)\| \le \frac{{e}^{\alpha T} - 1}{\alpha} \cdot \lim_{\Delta t \to 0} ( C_1 \Delta t + C_2 ) = \frac{e^{\alpha T} - 1}{\alpha} C_2.
\end{equation}
{
It demonstrates that while the discretization error vanishes as the sampling becomes instantaneous, the structural approximation error persists. The exponential term indicates that improving long-horizon predictive performance requires not only reducing local approximation errors, but also controlling the growth characteristics of the learned Koopman operators. In particular, excessive operator norms may amplify small modeling errors over time and degrade prediction accuracy. Motivated by this observation, the learning framework developed in the next section incorporates a stability-oriented regularization strategy that explicitly penalizes operator growth, thereby mitigating error amplification and improving long-term prediction performance.}


\section{Neural-Network-Based LTV Koopman Model}
\label{sec:NNKF}
{
The previous section established a control-oriented Koopman representation for nonlinear systems and developed an LTV approximation together with its associated error analysis. Building upon these results, this section presents a data-driven framework for constructing the proposed Koopman model directly from trajectory data. The framework aims to jointly learn three key components: the lifted-state representation, the time-varying Koopman dynamics, and the corresponding constraint representations. By learning these components within a unified architecture, the resulting model can capture complex nonlinear behaviors while preserving a linear structure in the lifted space that is suitable for prediction, optimization, and constraint handling. In addition, the error-propagation analysis developed in the previous section provides guidance for the design of a stability-oriented regularization strategy to improve long-horizon predictive performance. The overall architecture of the proposed neural-network-based Koopman framework is described in the following subsections.
}

\subsection{Lifting Function Parameterization and LTV Dynamics Learning}
To approximate the lifting function $\Psi$ in a data-driven manner, we parameterize it as
\begin{equation}
\Phi({\mathbf{x}}) = [\,1, {\mathbf{x}}^T, \psi_\theta({\mathbf{x}})^T\,]^T ,
\end{equation}
where the constant term allows the lifted dynamics to represent affine effects; the embedded original state ensures that the lifted space retains physical interpretability, allowing constraints and costs defined in the original state space to be imposed directly in the lifted domain; the learnable nonlinear feature block $\psi_\theta({x})$ is a multilayer perceptron whose depth, width, and activation functions can be adjusted, allowing the lifting network to flexibly capture the nonlinear features required by the Koopman representation. By augmenting the lifted state with both the raw physical variables and a rich set of learned nonlinear features, the lifting function enables the Koopman model to represent complex dynamics that are difficult to capture using hand-crafted basis functions. The resulting lifted dimension $n_z$ is determined automatically once the network is initialized.

Consistent with the time-varying linear Koopman model, we parameterize $A_k$ and $B_k$ as time-indexed linear maps acting on the lifted state. The learned lifted dynamics takes the form
\begin{equation} \label{eq:ltv_learn}
\hat {\mathbf{z}}_{k+1} = A_k \hat {\mathbf{z}}_k + B_k {\mathbf{u}}_k.
\end{equation}
Each operator pair $(A_k, B_k)$ is implemented as a bias-free linear layer whose parameters are learned independently for each time index. 
{
It is worth noting that although this structure introduces additional parameters scaling with the horizon length, it significantly expands the model's expressiveness while strictly preserving the linear structure required for downstream convex optimization. Crucially, this time-varying formulation still yields a global representation over the state space rather than relying on localizations around a single trajectory or specific initial conditions, which is further validated through experiments in Section~\ref{sec:simulation}. From a modeling perspective, the proposed framework can be interpreted as a generalization of conventional LTI Koopman formulations, where the time indices can be flexibly designed and the standard LTI representation emerges as a special case using a single time index.
}

\subsection{Constraint Lifting and Linear Constraint Representation}
In many control applications, the system is subject to nonlinear constraints of the form
\begin{equation}
\mathbf{h}(\mathbf{x}, \mathbf{u}) \le 0
\end{equation}
where $\mathbf{h}\in \mathbb{R}^{n+m} \to \mathbb{R}^p$, and it generally cannot be handled directly by linear or convex optimization methods. To enable constraint-aware control in the lifted space, we seek a representation in which these nonlinear constraints admit an approximately linear form with respect to the lifted variables.

Motivated by the Koopman framework, we treat the constraint functions as observables of the state–input pair and approximate them using linear operators acting on the lifted state. Specifically, at each time step $k$, we introduce a time-dependent linear constraint map
\begin{equation} \label{eq:cons_learn}
\hat {\mathbf{h}}_{k} = G_k \hat {\mathbf{z}}_k + H_k \mathbf{u}_k,
\end{equation}
where $G_k\in \mathbb{R}^{p\times n_z}$ and $H_k\in \mathbb{R}^{p\times m}$. This construction follows the same principle underlying Koopman lifting of dynamics: although the constraint functions are nonlinear in the original variables, they can often be approximated as linear functionals in a sufficiently rich lifted space. By sharing the neural-based lifted representation with the dynamics model, the learned observables are encouraged to support linear prediction of both the system evolution and the constraint values.

Allowing the constraint operators $G_k$ and $H_k$ to vary with time serves two purposes. First, it follows the LTV approximation used for the dynamics, ensuring consistency with the discretized bilinear Koopman structure and its local nature. Second, it accommodates large variations in constraint sensitivity along the trajectory, which would be difficult to capture using a single global linear map.

From an optimal control perspective, expressing constraints as linear functions of the lifted state and input is essential. It enables the resulting prediction and OCP to be formulated as a linear program, making it compatible with powerful linear programming solvers. Moreover, because the constraint approximation error is absorbed into the same lifted representation used for dynamics prediction, inaccuracies in constraint evaluation are naturally coupled to the state prediction error, rather than introducing an independent source of mismatch. For numerical stability and to avoid the introduction of artificial bias, each constraint operator pair $(G_k, H_k)$ is implemented as bias-free linear layers and learned independently for each time index. This ensures that constraint representations are consistent with the lifted dynamics.

Overall, the proposed constraint lifting strategy transforms nonlinear constraints into approximately linear expressions in the lifted coordinates while remaining consistent with the Koopman-based approximation of the dynamics. This unified lifted representation forms the basis for constraint-aware learning and control in the subsequent sections.

\subsection{Training Procedure and Loss Design}
Training proceeds by propagating the lifted dynamics over entire trajectories. Given state data $({\mathbf{x}}_0,\dots,{\mathbf{x}}_T)$, control inputs $({\mathbf{u}}_0,\dots,{\mathbf{u}}_{T-1})$, and constraint values $({\mathbf{h}}_0,\dots,{\mathbf{h}}_T)$, the model initializes the predicted lifted state as $\hat{{\mathbf{z}}}_0=\Phi({\mathbf{x}}_0)$ and propagated forward using the learned time-varying Koopman model~\eqref{eq:ltv_learn}.
Simultaneously, the predicted constraint values $\hat {\mathbf{h}}_k$ are obtained via equation~\eqref{eq:cons_learn}.
The loss function used to train the neural Koopman model is designed based on the multi-step error bounds derived in the previous sections. In particular, the analysis shows that long-horizon prediction accuracy depends critically on both the local approximation error and the stability of the time-varying Koopman operators. To reflect this dependence during learning, we adopt a multi-component training objective.

To improve both stability and long-horizon accuracy, we minimize a discounted multi-step loss for dynamics consistency
\begin{equation}
\mathcal{L}_{z}
= \frac{1}{M}\sum_{k=0}^{M-1} \gamma^k
\big\| \Phi(\mathbf{x}_k) - \hat {\mathbf{z}}_k \big\|^2,
\end{equation}
where $\gamma \in (0,1]$ is a discounted factor. Discounting later time steps mitigates the exponential amplification factor that appears in the theoretical error bound, while still enforcing consistency across the entire prediction horizon.

Similarly, to ensure that the lifted representation supports accurate constraint evaluation, we include a multi-step constraint consistency loss
\begin{equation}
\mathcal{L}_{h}
= \frac{1}{M}\sum_{k=0}^{M-1} \gamma^k
\big\|\mathbf{h}_k - \hat {\mathbf{h}}_k \big\|^2.
\end{equation}
Penalizing constraint mismatch along the full trajectory encourages the learned observables to encode constraint-relevant features in a manner that remains consistent over long horizons.

{
The multi-step error bound indicates that the induced norm of the Koopman operators governs the rate at which local approximation errors are amplified over the prediction horizon. To improve long-horizon prediction accuracy and mitigate excessive error propagation,} we include a regularization term on the deviation of $A_k$ from the identity:
\begin{equation} \label{eq:L_alpha}
\mathcal{L}_{\alpha}
= \frac{1}{M}\sum_{k=0}^{M-1}
\big\| A_k - I \big\|_F^2.
\end{equation}
{
This regularization penalizes overly aggressive state-transition operators and acts as a surrogate mechanism for controlling the growth factor appearing in the multi-step error bound. Consequently, it encourages stable long-horizon predictions while preserving sufficient flexibility to capture nonlinear dynamics.}

The full training objective is given by
\begin{equation}
\mathcal{L}
= \mathcal{L}_{z}
+ w_{h}\,\mathcal{L}_{h}
+ w_{\alpha}\,\mathcal{L}_{\alpha},
\end{equation}
where the weights $w_{h}$ and $w_{\alpha}$ are selected to balance prediction accuracy, constraint fidelity, and stability.

The model is trained end-to-end by minimizing the total loss $\mathcal{L}$ using stochastic gradient-based optimizer. Gradients are computed through the entire unrolled prediction horizon via backpropagation through time, allowing the lifting function and all time-varying Koopman and constraint operators to be optimized jointly. An adaptive learning-rate scheduler is employed to ensure stable convergence as the rollout horizon increases.

By explicitly designing the loss function to mitigate multi-step error accumulation, the proposed training scheme aligns the learned neural Koopman model with the underlying approximation theory. As a result, the learned lifted dynamics exhibit stability and accuracy over long horizons, making them well suited for downstream control and optimal control tasks.
{
It should be noted that because the lifting network $\Phi(\mathbf{x})$ and the time-varying operator sequences $(A_k, B_k)$ are optimized jointly under this multi-component objective, the resulting representation learning space is mathematically non-identifiable. Consequently, multiple equivalent combinations of latent coordinate transformations and operator sequences may yield identical input-output behavior. However, from a control-oriented perspective, unique identification of the true Koopman spectral components is not required. As long as the learned model minimizes $\mathcal{L}$ to ensure prediction and constraint accuracy, this non-uniqueness does not undermine the validity of the downstream convex optimization.
}

\section{Optimal Control with the Learned Model}
\label{sec:koopman_oc}
{
The previous sections established a neural-network-based LTV Koopman model capable of jointly learning lifted dynamics and constraint representations from data. This section leverages the learned Koopman model to construct a computationally tractable optimal control formulation for nonlinear constrained systems. Starting from a continuous-time nonlinear optimal control problem, we first reformulate the system dynamics and constraints in the learned Koopman coordinates using the LTV representation. The resulting lifted-space formulation preserves a linear prediction structure, allowing the original nonlinear optimal control problem to be approximated by a quadratic program (QP). Furthermore, the approximation-error analysis developed earlier is incorporated into the control design through a conservative constraint-tightening strategy, which explicitly accounts for modeling error when enforcing path constraints. The resulting framework provides a unified connection between Koopman-based learning, error-aware prediction, and constrained optimal control.
}

\subsection{Optimal Control in the Learned Lifted Space}

We consider the finite-horizon continuous-time optimal control problem
\begin{subequations} \label{eq:ocp}
\begin{align}
\min_{\mathbf{u}(\cdot)} \quad
& \int_0^T 
[
\mathbf{x}(t)^\top Q_c \mathbf{x}(t) + \mathbf{u}(t)^\top R_c \mathbf{u}(t)
]\,dt
+ \mathbf{x}(T)^\top Q_{f,c} \mathbf{x}(T)
\label{eq:oc_cost} \\
\text{s.t.} \quad
& \dot{\mathbf{x}}(t) = \mathbf{f}(\mathbf{x}(t)) + \sum_{i=1}^m \mathbf{g}_i(\mathbf{x}(t))\,\mathbf{u}_i(t),\quad \mathbf{x}(0)=\mathbf{x}_0
\label{eq:oc_dyn} \\
& \mathbf{h}(\mathbf{x}, \mathbf{u}) \le 0,\quad \forall t \in [0,T],
\label{eq:oc_con}
\end{align}
\end{subequations}
where $Q_c \succeq 0$, $R_c \succ 0$, and $Q_{f,c} \succeq 0$ are weighting matrices, the initial state $\mathbf{x}_0$ and the terminal time $T$ are given, and $\mathbf{h}(\mathbf{x}, \mathbf{u}) \le 0$ denotes general nonlinear constraints. Problem \eqref{eq:ocp} is nonconvex in general due to nonlinear dynamics and constraints, making it challenging to solve.

Using the learned lifting function $\Phi$, the initial system state is transferred to the lifted space as $\hat{\mathbf{z}}_0 = \Phi(\mathbf{x}_0)$.
Under piecewise-constant control inputs and sampling period $\Delta t$, the lifted dynamics is approximated by the linear time-varying Koopman model as
\begin{equation}
\hat{\mathbf{z}}_{k+1} = A_k \hat{\mathbf{z}}_k + B_k \mathbf{u}_k,
\label{eq:lifted_dyn}
\end{equation}
where the operators $(A_k,B_k)$ are learned using the neural Koopman framework described in the previous section.

Consistent with the continuous-time cost \eqref{eq:oc_cost}, we define the discrete-time lifted quadratic cost
\begin{equation} \label{eq:lifted_cost}
\ell_k(\hat{\mathbf{z}}_k,\mathbf{u}_k) = \hat{\mathbf{z}}_k^\top Q \hat{\mathbf{z}}_k + \mathbf{u}_k^\top R \mathbf{u}_k,
\qquad
\ell_f(\hat{\mathbf{z}}_N) = \hat{\mathbf{z}}_N^\top Q_f \hat{\mathbf{z}}_N,
\end{equation}
where the weighting matrices $Q,R,Q_f$ are constructed following the standard discretization approach.
Moreover, the nonlinear constraints $\mathbf{h}(\mathbf{x},\mathbf{u})\le \mathbf{0}$ are approximated in the lifted space by linear constraints
\begin{equation}
G_k \hat{\mathbf{z}}_k + H_k \mathbf{u}_k \le 0,
\label{eq:lifted_con}
\end{equation}
as introduced in the neural Koopman framework.

\subsection{Constraint-Tightened Convex Programming}

The lifted dynamics \eqref{eq:lifted_dyn} and constraints \eqref{eq:lifted_con} are linear with respect to the lifted state and control. However, since the Koopman model provides only an approximation of the true lifted dynamics, constraint satisfaction for the predicted lifted state does not automatically guarantee satisfaction for the true system.
To address this issue, we incorporate a conservative constraint tightening based on the Koopman prediction error bounds.

From the multi-step error propagation analysis, the prediction error $\mathbf{e}_k = \mathbf{z}_k - \hat{\mathbf{z}}_k$ admits an upper bound that depends on the accumulated single-step approximation errors and the growth properties of the Koopman operators, denoted as $\eta_k$. {
While the theoretical analysis establishes the structure of the error propagation bound, the magnitude of the single-step approximation error is an a posteriori property of the learned Koopman model and therefore cannot be determined independently of model validation. Consequently, the single-step approximation error is estimated from validation prediction residuals and then propagated using the multi-step error model derived in Section~\ref{sec:error propagation} to obtain an estimation $\hat\eta_k$.}

Substituting $\mathbf{z}_k = \hat{\mathbf{z}}_k + \mathbf{e}_k$ into the lifted constraint \eqref{eq:lifted_con} yields
\[
G_k \hat{\mathbf{z}}_k + H_k \mathbf{u}_k + G_k \mathbf{e}_k \le \mathbf{0}.
\]
Using the Cauchy--Schwarz inequality, the constraint perturbation induced by the prediction error can be bounded as
\[
G_k \mathbf{e}_k
\le
\|G_k\|\,\|\mathbf{e}_k\|
\le
\|G_k\|\,\hat{\eta}_k.
\]
Therefore, a conservative strategy for improving constraint satisfaction is to enforce the tightened constraint
\begin{equation}
G_k \hat{\mathbf{z}}_k + H_k \mathbf{u}_k
\le
-\|G_k\|\,\hat{\eta}_k.
\label{eq:tight_con}
\end{equation}
{ Replacing \eqref{eq:lifted_con} with \eqref{eq:tight_con} yields a conservatively tightened quadratic program that explicitly accounts for Koopman model prediction errors during constraint enforcement. The resulting control solution provides an additional margin against constraint violations and improves the reliability of constraint satisfaction in the presence of Koopman model uncertainty.
}

Using the tightened constraint \eqref{eq:tight_con}, lifted dynamics \eqref{eq:lifted_dyn}, and cost \eqref{eq:lifted_cost}, the entire optimal control problem can be reformulated as
\begin{subequations} \label{eq:qp}
\begin{align}
\min_{\{\mathbf{u}_k\}_{k=0}^{N-1}} \quad
& \sum_{k=0}^{N-1}
\left(
\mathbf{z}_k^\top Q \mathbf{z}_k + \mathbf{u}_k^\top R \mathbf{u}_k
\right)
+ \mathbf{z}_N^\top Q_f \mathbf{z}_N
\label{eq:qp_cost} \\
\text{s.t.} \quad
& \mathbf{z}_{k+1} = A_k \mathbf{z}_k + B_k \mathbf{u}_k, \quad k=0,\dots,N-1,
\label{eq:qp_dyn} \\
& \mathbf{z}_0 = \Phi(\mathbf{x}_0),
\label{eq:qp_init} \\
& G_k \mathbf{z}_k + H_k \mathbf{u}_k \le  - \|G_k\| \, \hat\eta_k, \quad k=0,\dots,N-1.
\label{eq:qp_con}
\end{align}
\end{subequations}
Due to the quadratic cost and linear dynamics and constraints, problem \eqref{eq:qp} is a convex quadratic program.

The proposed Koopman-based optimal control formulation offers several advantages. First, by leveraging the learned LTV Koopman representation, the original nonlinear optimal control problem can be approximated by a convex quadratic program, enabling computationally efficient prediction and control. Second, the unified lifted-space framework allows system dynamics and constraint representations to be learned and utilized within a consistent optimization formulation. Third, the proposed constraint-tightening strategy provides an error-aware mechanism for constraint enforcement. Overall, the framework establishes a systematic connection between data-driven Koopman learning and constrained optimal control, combining modeling flexibility, computational tractability, and explicit consideration of approximation uncertainty.

\section{Simulation Results}
\label{sec:simulation}
{
This section evaluates the performance of the proposed framework using two representative aerospace optimal control problems. We first consider a low-thrust maximum-radius orbit transfer problem~\cite{Ross2007LowThrustHT}, which is characterized by nonlinear control-affine dynamics with a nonlinear equality-constrained target set. This problem serves as a baseline to systematically evaluate the approximation of nonlinear dynamics and constraints across model configurations, directly illustrating the impact of model accuracy on the optimal control solutions. To further assess the applicability and robustness of the proposed framework in a more challenging setting, we then consider a six-degree-of-freedom (6-DOF) atmospheric entry guidance problem~\cite{Wan2022F, Pei2023A}. This problem exhibits nonlinearly coupled translation-rotation dynamics and multiple nonlinear inequality path constraints, providing a rigorous scenario for evaluating the framework's adaptability and the effectiveness of the constraint-tightening strategy. For both problems, the prediction performance of different Koopman model variants is analyzed through ablation studies to isolate the contribution of each component. Subsequently, the Koopman-based optimal control formulation across different variants is evaluated and benchmarked against standard NLP solvers, highlighting the advantages of the proposed framework in terms of solution accuracy, robustness, and computational efficiency. All experiments were conducted within a unified Python environment on a standard laptop CPU.

\subsection{Maximum Radius Orbit Transfer Problem}
We first evaluate the proposed framework on the maximum radius orbit transfer problem~\cite{Ross2007LowThrustHT}. The dynamics of the low-thrust orbit transfer are given by
\begin{equation}\label{eq:orbit_dyn}
\begin{aligned}
&\dot{r}=v_r,\\
&\dot{v}_r=\frac{v_\theta^2}{r} - \frac{\mu}{r^2} + u_r,\\
&\dot{v}_\theta=-\frac{v_r v_\theta}{r} + u_\theta,
\end{aligned}
\end{equation}
where $r, v_r, v_\theta$ represent the radial position, radial velocity and transverse velocity. Control inputs $u_r=T_r/m$ and $u_t=T_t/m$ correspond to the radial and transverse components of the thrust. Canonical units are employed to scale the orbit dynamics, where $r$ is expressed in astronomical unit (AU), $v_r$ and $v_\theta$ are in AU per time unit (TU), and the canonical gravitational constant is normalized to $\mu=1$. The normalized control vector (non-dimensional acceleration) is bounded within $\mathcal{U}=[u_{\min}, u_{\max}]^2$. The terminal constraints at $t_f$ ensure that the spacecraft injects into a circular orbit, which is characterized by the target set
\begin{equation}\label{eq:orbit_con}
\mathcal{X}_f = \left\{ [r, v_r, v_\theta]^\top \in \mathbb{R}^3 \;\middle|\; v_r = 0, \; rv_\theta^2-\mu = 0 \right\}.
\end{equation}

The maximum radius orbit transfer problem aims to maximize the spacecraft's radial distance using continuous thrust, which is formulated as the following optimal control problem:
\begin{subequations} \label{eq:entry}
\begin{align}
\min_{\mathbf{u}(t)}\quad& -r(t_f)\\
\text{subject to}\quad& \dot{\mathbf{x}} = \mathbf{f}(\mathbf{x}, \mathbf{u}),\\
&\mathbf{x}(t_0) = \mathbf{x}_0,\quad \mathbf{x}(t_f) \in \mathcal{X}_f,\\
&\mathbf{u}(t) \in \mathcal{U},\quad t\in[t_0,\ t_f],
\end{align}
\end{subequations}
where $\mathbf{x} = [r,\ v_r,\ v_\theta]^T$ is the state vector, $\mathbf{u} = [u_r,\ u_\theta]^T$ is the control vector, $\mathbf{f}$ denotes the dynamics in \eqref{eq:orbit_dyn}.

\subsubsection{Prediction Performance Evaluation}
To systematically evaluate the contributions of each component, we conduct an ablation study on the maximum-radius orbit transfer problem. Four models are trained to learn the control-affine dynamics in~\eqref{eq:orbit_dyn} together with the nonlinear circular-orbit constraint in~\eqref{eq:orbit_con}. Specifically, we consider (i) the proposed LTV Koopman model with a neural lifting function (LTV-NN), (ii) an LTI Koopman model with a neural lifting function (LTI-NN), (iii) an LTV Koopman model with a fixed polynomial lifting function (LTV-Poly), and (iv) an LTI Koopman model with a fixed polynomial lifting function (LTI-Poly). The LTI-Poly and LTI-NN baselines adopt representative LTI Koopman predictor structures with fixed and neural-network-based lifting functions, as studied by Korda and Mezi{\'c}~\cite{Korda2018L} and Han et al.~\cite{han2020deep}, respectively. For the fixed lifting, polynomials up to fourth order are included in the observable dictionary. The neural lifting is implemented using a three-layer multilayer perceptron (MLP) with the same output dimension as the polynomial lifting. Training data are generated by randomly sampling states from the domain $r\in[0.5,1.5]$, $v_r\in[-0.5,0.5]$, $v_\theta\in[0.5,1.5]$, and propagating the dynamics for $1.5$ TU with a uniform sampling interval of $\Delta t=0.05$ TU under random controls satisfying $u_r,u_\theta \in[-0.1,0.1]$. A total of $5000$ trajectories are generated and split into training, validation, and test sets with ratios of $80\%$, $10\%$, and $10\%$, respectively. For a fair comparison across different variants, the operator matrices are uniformly initialized such that the state transition matrix ($A_k$ or $A$) is a identity matrix and the control matrix ($B_k$ or $B$) is a zero matrix. All models are trained using the Adam optimizer~\cite{kingma2015adam} with identical optimization settings.

To evaluate the prediction accuracy, the trained models are propagated over the entire horizon using initial states and control inputs from the test set, and then compared with the ground truth trajectories. The resulting relative root mean square error (RMSE) for dynamics and constraint predictions is shown in Fig.~\ref{fig:orbit_pred}. The proposed LTV Koopman model with neural-network-based lifting achieves the lowest prediction error throughout the horizon, with both dynamics and constraint errors remaining below $2\%$ at the final time. For both LTV and LTI formulations, neural lifting consistently outperforms the fixed polynomial lifting, indicating that the learned observables provide a more expressive representation of the nonlinear dynamics and constraints. By jointly optimizing the lifting function, Koopman dynamics, and constraint representations, the neural lifting adapts the lifted space to the prediction and constraint-linearization objectives, leading to improved accuracy. In addition, the LTV models exhibit substantially slower error accumulation than their LTI counterparts, demonstrating the advantage of time-varying operators in capturing the evolving nonlinear dynamics over long horizons. Overall, the results confirm that both the neural lifting function and the LTV formulation contribute to improved long-horizon prediction accuracy, with their combination (LTV-NN) yielding the best performance.

\begin{figure}[htbp]
    \centering
    \includegraphics[width=0.8\linewidth]{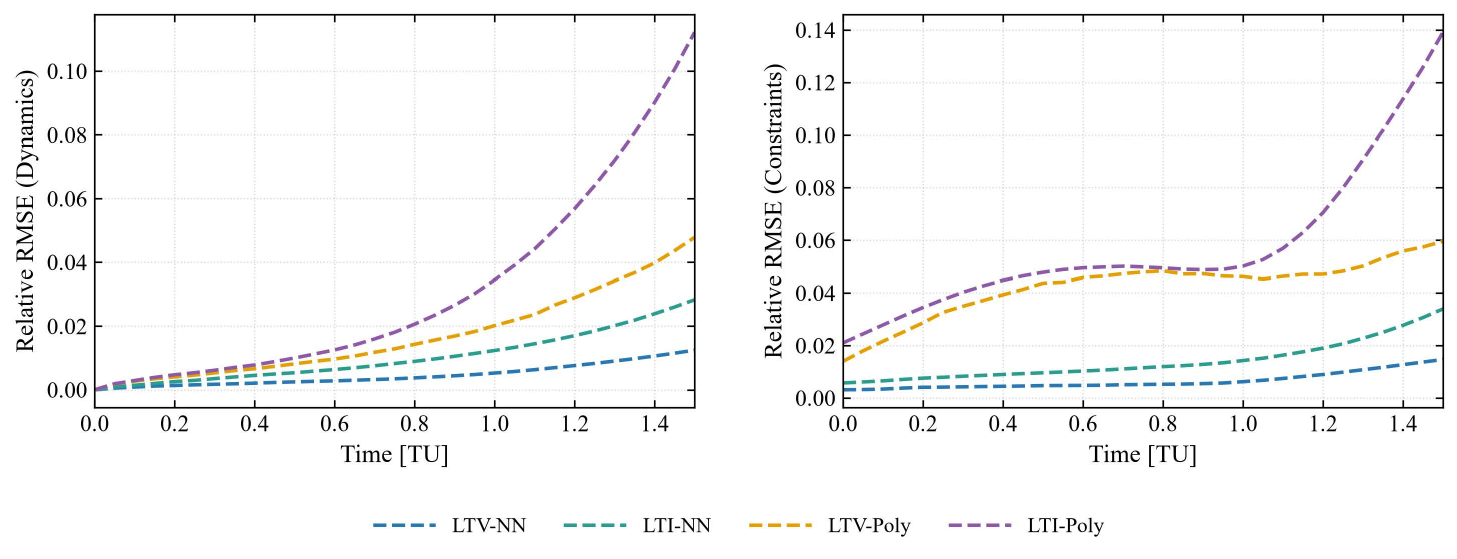}
    \caption{Prediction errors across different models for the orbit transfer dynamics}
    \label{fig:orbit_pred}
\end{figure}

The effect of the stability regularization is further investigated for both the LTV and LTI Koopman models with neural lifting. The stability regularization term~\eqref{eq:L_alpha} is introduced with different weights $w_\alpha$ while keeping all other training settings unchanged. The resulting prediction errors are shown in Fig.~\ref{fig:orbit_weight}. For the LTV model, introducing stability regularization slightly increases the short-horizon prediction error, indicating a reduced flexibility in fitting the dynamics. However, moderate values of $w_\alpha$ significantly suppress the exponential error amplification for long-horizon without substantially compromising local prediction accuracy. For the LTI model, increasing $w_\alpha$ consistently degrades short-horizon prediction accuracy, while providing very limited improvement in long-horizon performance. A possible explanation is that the time-invariant Koopman operator already imposes a restrictive approximation of the nonlinear dynamics, and the additional regularization further constrains the model capacity, causing the increase in approximation error to outweigh the benefits of improved operator stability.

\begin{figure}[htbp]
    \centering
    \includegraphics[width=0.8\linewidth]{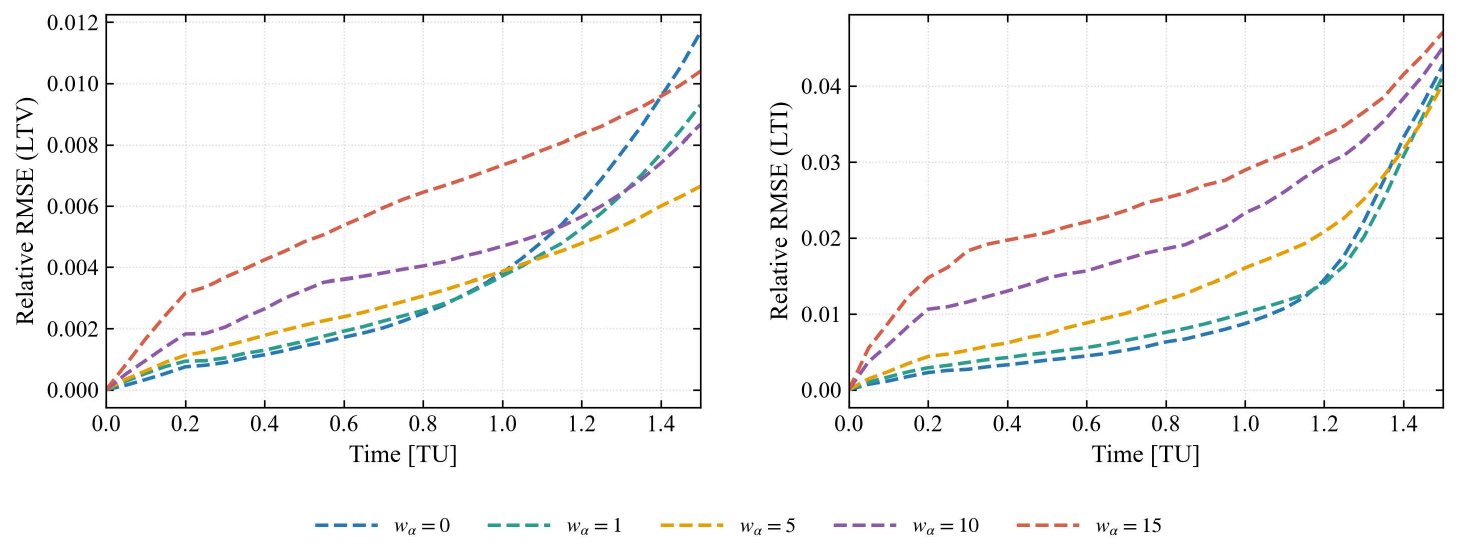}
    \caption{Prediction errors with different $w_\alpha$ for the orbit transfer dynamics}
    \label{fig:orbit_weight}
\end{figure}

\subsubsection{Optimal Control Evaluation}
To further evaluate the proposed optimal control framework, the four Koopman model variants considered in the prediction study are applied to the maximum-radius orbit transfer problem. For each model, both the nonlinear dynamics and the circular-orbit terminal constraint are transformed into linear representations in the lifted space, yielding a constrained linear programming problem. Since the terminal constraint set $\mathcal{X}_f$ is defined as an equality level set, it is enforced directly in the lifted optimization problem without constraint tightening, directly reflecting the impact of constraint accuracy on optimal solutions. To benchmark performance, we compare the Koopman-based approaches against two standard NLP solvers based on the interior-point method (NLP-IP) and the sequential quadratic programming method (NLP-SQP). 

The results are summarized in Table~\ref{tab:orbit}, where each method is evaluated over 50 Monte Carlo cases with initial states randomly sampled from the training domain and randomized initial guesses for the NLP solvers. Representative solutions are shown in Fig.~\ref{fig:orbit_sol}. The results show that, for the Koopman-based approaches, control performance is strongly correlated with model accuracy. More accurate lifted representations yield solutions that are closer to the NLP reference, exhibit smaller terminal constraint violations, and achieve higher success rates. This trend is evident across the four Koopman variants, with LTV-NN achieving the best overall performance. As model accuracy decreases from LTV-NN to LTI-Poly, objective optimality and constraint satisfaction all degrade accordingly. In some cases, inaccurate lifted constraint distort the terminal feasible set, ultimately causing optimization failure. In contrast, although the NLP methods achieve higher solution accuracy, the occasional convergence failures reflect the sensitivity of nonlinear optimization to initial guess and local convergence behavior. This sensitivity also results in substantially larger computational costs and variability in solution time. Overall, the results demonstrate that improved lifted-space representations in the proposed framework lead directly to higher-quality optimal control solutions, narrowing the performance gap with nonlinear optimization while retaining the robustness and computational efficiency of Koopman-based formulations.

\begin{figure}[htbp]
    \centering
    \includegraphics[width=\linewidth]{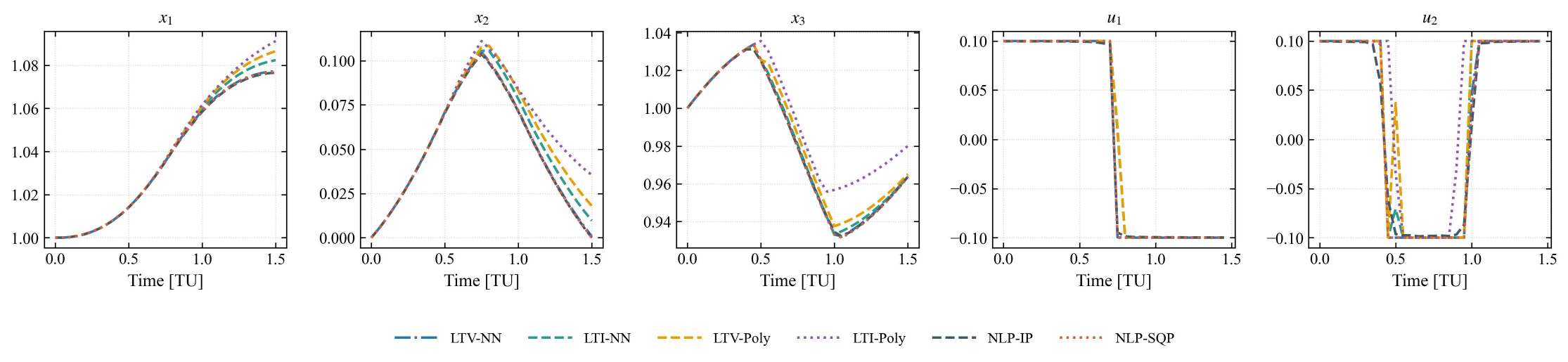}
    \caption{Representative solutions comparison for the maximum radius orbit transfer problem}
    \label{fig:orbit_sol}
\end{figure}

\begin{table}[htbp]
\centering
\caption{Performance metrics for the maximum radius orbit transfer problem based on Monte Carlo simulations}
\label{tab:orbit}
\begin{tabular}{lcccccc}
\hline
\textbf{Method} & \textbf{Avg. Relative Objective} & \textbf{Avg. Constraint Violation} & \textbf{Avg. Time (s)} & \textbf{Success Rate (\%)} \\ \hline
NLP-IP    & 1 (Ref.)    & 4.5e-8     & 9.27    & 98\\
NLP-SQP    & 1.0006     & 8.7e-8    & 2.19    & 92\\
LTV-NN    & 1.0057    & 2.1e-3    & 0.097    & 100\\
LTI-NN    & 1.0087    & 8.9e-3    & 0.090    & 100\\ 
LTV-Poly    & 1.0093    & 9.1e-3    & 0.101    & 96\\ 
LTI-Poly    & 1.0142    & 4.5e-2    & 0.103    & 82\\ \hline
\end{tabular}
\end{table}
}

\subsection{6-DOF Entry Guidance Problem}


The proposed framework is next evaluated in the 6-DOF atmospheric entry guidance problem~\cite{Wan2022F, Pei2023A}. The atmospheric entry vehicle is modeled as a rigid body subject to aerodynamic and gravitational forces. The translational motion is represented by the radial distance $r$, longitude $\lambda$, latitude $\phi$, velocity magnitude $V$, flight-path angle $\gamma$, and heading angle $\psi$. The rotational motion is parameterized by the Euler angles $(\theta_1,\theta_2,\theta_3)$ and body angular rates $(\omega_x,\omega_y,\omega_z)$. The control inputs consist of the angle of attack $\alpha$ and side-slip angle $\beta$, which determine the orientation of the aerodynamic force vector and influence both the translational and rotational dynamics.

The translational equations of motion follow the standard spherical-coordinate formulation projected along the velocity-aligned reference frame:
\begin{subequations}\label{eq:eom_trans}
\begin{align}
\dot{r} &= V \sin\gamma, \quad \dot{\phi} = \frac{V\cos\gamma\sin\psi}{r}, \quad \dot{\lambda} = \frac{V\cos\gamma\cos\psi}{r\cos\phi}, \\
\dot{V} &= -D - \frac{\mu\sin\gamma}{r^2}, \quad \dot{\gamma} = \frac{L\cos\sigma}{V} + \frac{V\cos\gamma}{r} - \frac{\mu\cos\gamma}{r^2\,V}, \\
\dot{\psi} &= \frac{L\sin\sigma}{V\cos\gamma} + \frac{V\cos\gamma\cos\psi \tan\phi}{r},
\end{align}
\end{subequations}
where $\sigma$ is the bank angle. The aerodynamic lift $L$ and drag $D$ are modeled as:
\begin{equation}
L = \frac{1}{2} \rho(r) V^{2} S C_L(\alpha), \quad D = \frac{1}{2} \rho(r) V^{2} S C_D(\alpha),
\end{equation}
with $S$ being the reference area, $\rho(r)$ the atmospheric density, and the aerodynamic coefficients given by $C_L = C_{L_0} + C_{L_1}\alpha$ and $C_D = C_{D_0} + C_{D_1}\alpha + C_{D_2}\alpha^2$.

The attitude kinematics and rigid-body rotational dynamics are parameterized using Euler angles $(\theta_1,\theta_2,\theta_3)$ and body angular rates $(\omega_x,\omega_y,\omega_z)$:
\begin{subequations}\label{eq:attitude_dynamics}
\begin{align}
\begin{bmatrix} \dot{\theta}_{1} \\ \dot{\theta}_{2} \\ \dot{\theta}_{3} \end{bmatrix} &= \frac{1}{\cos\theta_2} 
\begin{bmatrix} 
-\omega_x\sin\theta_3 + \omega_z\cos\theta_3 \\ 
\omega_x\cos\theta_3\cos\theta_2 + \omega_z\sin\theta_3\cos\theta_2 \\ 
\omega_x\sin\theta_3\sin\theta_2 + \omega_y\cos\theta_2 - \omega_z\cos\theta_3\sin\theta_2 
\end{bmatrix}, \label{eq:eom_euler} \\
\begin{bmatrix} \dot{\omega}_x \\ \dot{\omega}_y \\ \dot{\omega}_z \end{bmatrix} &= 
\begin{bmatrix} 
(M_x - (J_2 - J_3) \omega_y \omega_z)/J_1 \\ 
(M_y - (J_3 - J_1) \omega_x \omega_z)/J_2 \\ 
(M_z - (J_1 - J_2) \omega_x \omega_y)/J_3 
\end{bmatrix}, \label{eq:eom_rates}
\end{align}
\end{subequations}
where $J_i$ ($i=1,2,3$) are the principal moments of inertia. The aerodynamic moment vector is $\mathbf{M}_a = [M_x, M_y, M_z]^\top = \mathbf{r}_{ac} \times \mathbf{F}_a$, where $\mathbf{r}_{ac}$ is the position vector of the aerodynamic center, and $\mathbf{F}_a = [W_x, W_y, W_z]^\top$ represents the body-axis aerodynamic forces resolved as:
\begin{equation}\label{eq:ld_force}
\begin{aligned}
W_x &= \sin{\beta} (D \cos{\alpha} + L \sin{\alpha}), \quad W_y = -\cos{\beta} (D \cos{\alpha} + L \sin{\alpha}), \quad W_z = L \cos{\alpha} - D \sin{\alpha}.
\end{aligned}
\end{equation}
All variables are nondimensionalized using the Mars radius $R_0=3389.5\text{ km}$, gravitational acceleration $g_0=3.711\text{ m/s}^2$, and characteristic velocity $\sqrt{R_0g_0}$.

The vehicle must satisfy aerothermal and structural constraints during atmospheric flight. The path constraints on the stagnation-point heating rate $\dot{Q}$ and structural normal load are given by
\begin{equation} \label{eq:path_cons}
\dot Q = k_Q \sqrt{\frac{\rho}{R_{\text{nose}}}} V^{3.15} \le \dot Q_{\max}, \quad \sqrt{L^2 + D^2} \le n_{\max},
\end{equation}
where $k_Q$ is the atmospheric constant, $R_{\text{nose}}$ is the nose radius, and $\dot Q_{\max}$ and $n_{\max}$ are the maximum allowable limits. The control variables (angle of attack $\alpha$ and side-slip angle $\beta$) are bounded within the admissible control set $\mathcal{U}$, and the terminal altitude requirement restricts the final state to the target set $\mathcal{X}_f$, defined respectively as
\begin{equation} \label{eq:entry_sets}
\begin{aligned}
\mathcal{U} &= \left\{ \mathbf{u} \in \mathbb{R}^2 \;\middle|\; \alpha^L \le \alpha \le \alpha^U, \; \beta^L \le \beta \le \beta^U \right\}, \\
\mathcal{X}_f &= \left\{ \mathbf{x}(t_f) \in \mathbb{R}^{12} \;\middle|\; r_f^L \le r(t_f) \le r_f^U \right\},
\end{aligned}
\end{equation}
where $\mathbf{x} = [r, \lambda, \phi, V, \gamma, \psi, \theta_1, \theta_2, \theta_3, \omega_x, \omega_y, \omega_z]^T$ is the state vector, $\mathbf{u} = [\alpha, \beta]^T$ is the control vector, the superscripts $L$ and $U$ denote the lower and upper bounds of the respective variables.

To reduce the propellant consumption required in the subsequent powered-descent phase, the final velocity at the fixed terminal time is minimized. Therefore, the entry optimal control problem is formulated as:
\begin{subequations} \label{eq:entry}
\begin{align}
\min_{\mathbf{u}(t)}\quad& V(t_f)\\
\text{subject to}\quad& \dot{\mathbf{x}} = \mathbf{f}(\mathbf{x}, \mathbf{u}),\label{eq:entry_f}\\
&\mathbf{h}\big(\mathbf{x}(t),\mathbf{u}(t)\big)\le 0,\label{eq:entry_h}\\
&\mathbf{x}(t_0) = \mathbf{x}_0,\quad \mathbf{x}(t_f) \in \mathcal{X}_f,\\
&\mathbf{u}(t) \in \mathcal{U},\quad t\in[t_0,\ t_f],
\end{align}
\end{subequations}
where $\mathbf{f}$ collects the translational and rotational dynamics \eqref{eq:eom_trans}--\eqref{eq:ld_force}, $\mathbf{h}$ includes the path constraints \eqref{eq:path_cons}. { Although the entry guidance problem is not a strictly control-affine system, this challenging application is used to validate the data-driven adaptability of the
proposed framework beyond the idealized control-affine setting, as well as the effectiveness of the
constraint-tightening strategy under stronger nonlinearity and more complex path constraints.}

{
\subsubsection{Prediction Performance Evaluation}
To evaluate the prediction performance of the proposed framework on the 6-DOF entry guidance problem, we conduct a systematic comparison across the same four Koopman model variants under a unified experimental setup. All models are trained to learn the lifted representation of the nonlinear entry dynamics~\eqref{eq:entry_f} together with associated path constraints~\eqref{eq:entry_h}, following the same architecture and training protocol as in the orbit transfer case.} For the entry problem, a dataset of $5000$ trajectories is generated under randomized initial conditions in ranges: $r_0 \in [98, 102]$ km, $\lambda_0 \in [110^\circ, 130^\circ]$, $\phi_0 \in [30^\circ, 50^\circ]$, $V_0 \in [4600, 4800]$ m/s, $\gamma_0 \in [0^\circ, 20^\circ]$, $\psi_0 \in [-10^\circ, 10^\circ]$, and control profiles within the admissible sets: $\alpha, \beta \in [-20^\circ, 20^\circ]$. Each trajectory is simulated over $T=20$ steps with a uniform sampling interval of $\Delta t=5$ s. {The dataset is split into training, validation, and test sets with ratios of $80\%$, $10\%$, and $10\%$, respectively. All models are trained using the Adam optimizer with identical hyperparameters as in the orbit transfer case.

To evaluate the prediction accuracy, the trained models are propagated using initial states and control inputs from the test set, and then compared with the ground truth trajectories. The resulting relative RMSE for both dynamics and constraint predictions is shown in Fig.~\ref{fig:entry_pred}. Overall, the results are consistent with those observed in the orbit transfer problem. The proposed LTV Koopman model with neural lifting achieves the lowest prediction error across the entire prediction horizon, while neural lifting consistently outperforms the fixed polynomial dictionary for both LTV and LTI formulations. Compared with the orbit transfer case, the performance gap between the models becomes more pronounced. In the 6-DOF entry problem, the system exhibits stronger nonlinear couplings, rapidly varying aerodynamic effects, and more complicated nonlinear constraints, which are less effectively captured by time-invariant operators and hand-crafted lifting functions. As a result, the LTV formulation with neural networks (LTV-NN) provides a more accurate prediction of dynamics and constraints along the trajectory, particularly in more complex nonlinear systems.

\begin{figure}[htbp]
    \centering
    \includegraphics[width=0.8\linewidth]{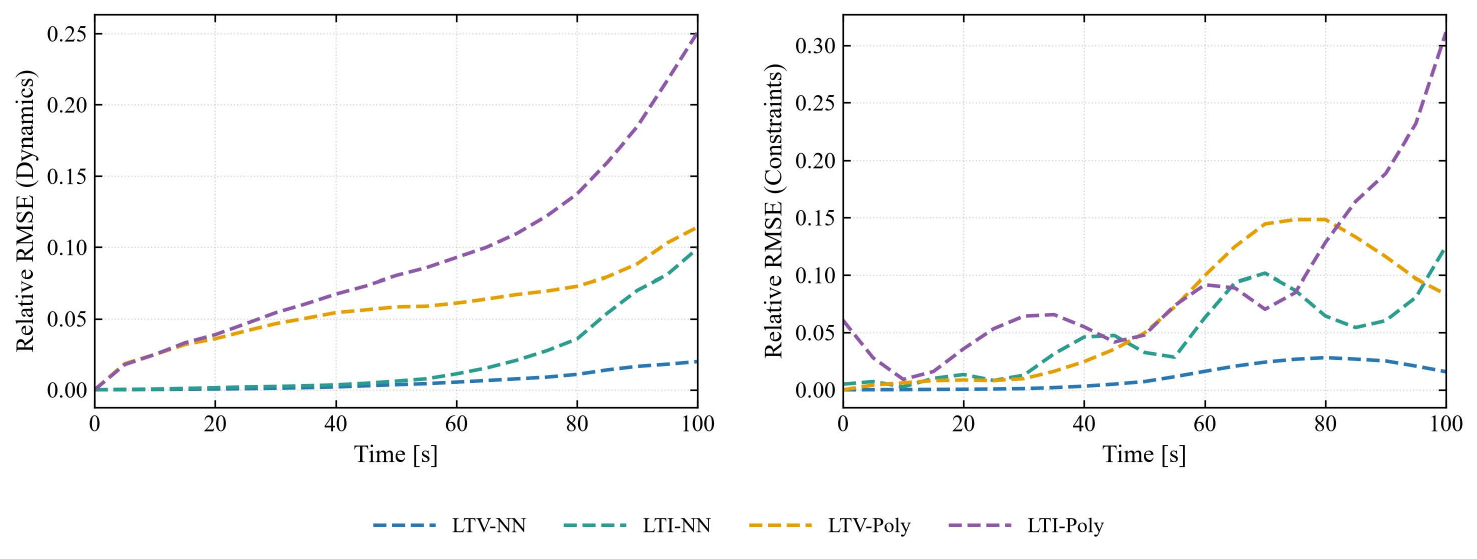}
    \caption{Prediction errors across different models for the entry guidance problem}
    \label{fig:entry_pred}
\end{figure}

We further investigate the effect of stability regularization for both LTV and LTI Koopman models with neural lifting on the 6-DOF entry guidance problem. The same regularization formulation is applied with varying weights $w_\alpha$, while keeping all other training settings unchanged. The results in Fig.~\ref{fig:entry_weight} are consistent with those observed in the orbit transfer problem. For the LTV model, the effect of stability regularization follows a similar trend as in the orbit transfer case, reflecting a trade-off between short-horizon accuracy and long-horizon stability. For the LTI model, its impact remains limited in improving long-horizon performance, while introducing additional degradation in short-horizon prediction accuracy as the regularization strength increases. This suggests that the expressiveness limitation of time-invariant operators remains the dominant factor in the entry problem, which restricts the effectiveness of additional stability constraints. Overall, the results indicate that stability regularization consistently improves long-horizon prediction performance under the LTV formulation in the entry guidance problem.

\begin{figure}[htbp]
    \centering
    \includegraphics[width=0.8\linewidth]{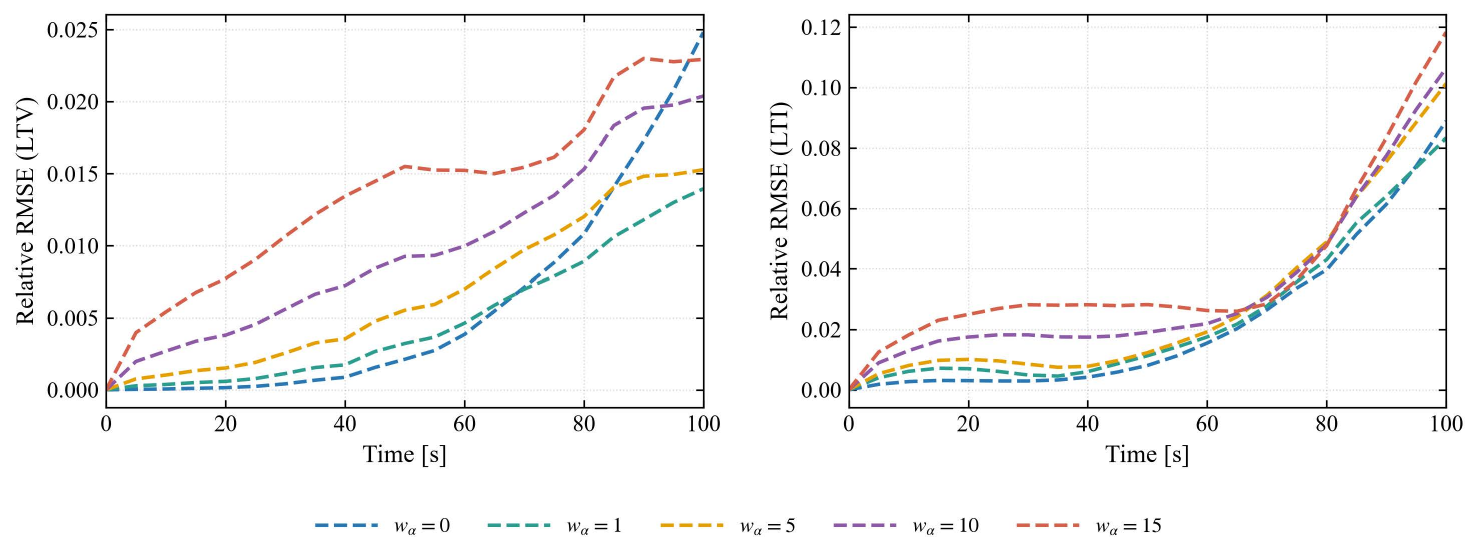}
    \caption{Prediction errors with different $w_\alpha$ for the entry guidance problem}
    \label{fig:entry_weight}
\end{figure}

\subsubsection{Optimal Control Evaluation}
To assess the optimal control performance, the learned Koopman models are applied to the entry guidance problem and compared with NLP methods. In addition to the four Koopman model variants considered in the prediction study, we further evaluate the effect of the proposed constraint tightening strategy by comparing with the non-tightened LTV Koopman model (LTV-NT). Each method is evaluated over 50 Monte Carlo cases with initial states randomly sampled from the training domain and randomized initial guesses for the NLP solvers. The results are summarized in Table~\ref{tab:entry}, with representative solutions shown in Fig.~\ref{fig:entry_sol}.

Among all Koopman-based approaches, LTV-NN achieves the best overall performance, yielding solutions closest to the NLP reference with a much higher success  than NLP. 
In comparison, LTI-NN exhibits noticeably larger objective degradation, constraint gap, and a substantially lower success rate. The two polynomial lifting models fail to generate feasible solutions across all Monte Carlo trials. Compared with the orbit transfer problem, the performance gap between different Koopman variants becomes significantly larger, indicating that the higher-dimensional and strongly nonlinear entry dynamics place much stricter requirements on prediction accuracy. In particular, prediction errors not only degrade the precision of lifted dynamics approximation but also enlarge the tightening margins required for constraint satisfaction, leading to increasingly conservative optimization results and even optimization failure for insufficiently accurate models.

The influence of constraint tightening can be observed by comparing LTV-NN and LTV-NT. Without tightening, the optimizer achieves a slightly better objective value than the NLP reference on average. However, this is caused by model approximation errors that make the lifted optimization problem overly optimistic, resulting in constraint violations when the solutions are evaluated using the original dynamics. As illustrated in Fig.~\ref{fig:entry_sol}, the non-tightened solution exceeds the heating-rate limits, whereas the tightened formulation maintains a safe margin throughout the trajectory. Although constraint tightening leads to a slightly lower success rate and introduces a small loss in objective optimality, it effectively compensates for prediction and constraint approximation errors, thereby improving the reliability of the obtained solutions.

Since the entry guidance problem is significantly more challenging due to its higher dimensionality and stronger nonlinearities, it make the optimization landscape considerably more difficult. As a result, the performance of NLP methods is more sensitive to the initial guess quality and local convergence behavior, as reflected in the relatively lower success rates and the increased computational time across Monte Carlo trials. In contrast, the proposed framework reformulates the optimal control problem as a linear program in the lifted space, enabling substantially more efficient computation and eliminating the dependence on carefully designed initial guesses. Together with the constraint-tightening strategy, this allows the method to maintain reliable feasibility while preserving competitive solution quality, thereby achieving a balance between optimality, efficiency, and robustness for complex nonlinear optimal control.

\begin{figure}[htbp]
    \centering
    \includegraphics[width=\linewidth]{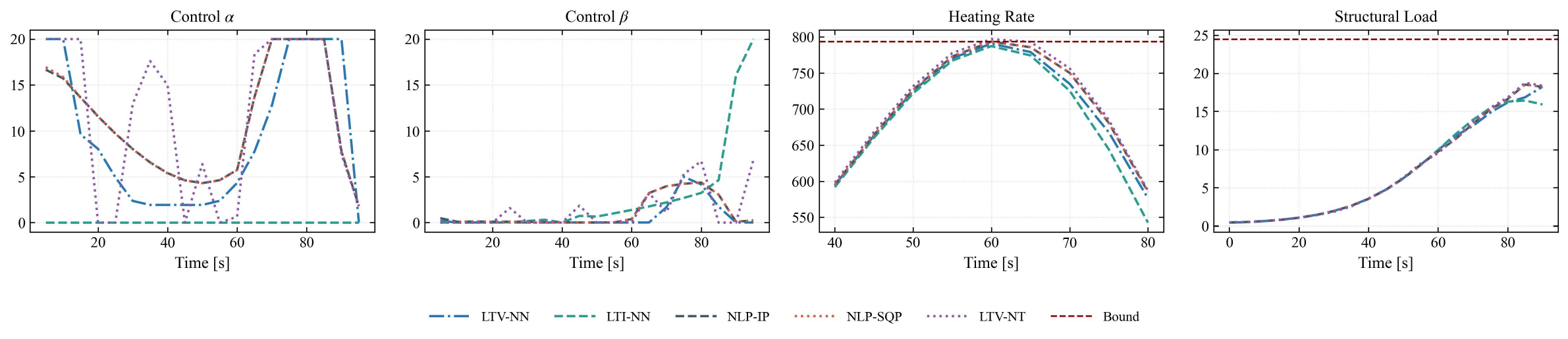}
    \caption{Representative solutions comparison for the entry guidance problem}
    \label{fig:entry_sol}
\end{figure}

\begin{table}[htbp]
\centering
\caption{Performance metrics for the entry guidance problem based on Monte Carlo simulations}
\label{tab:entry}
\begin{tabular}{lcccccc}
\hline
\textbf{Method} & \textbf{Avg. Relative Objective} & \textbf{Avg. Constraint Gap} & \textbf{Avg. Time (s)} & \textbf{Success Rate (\%)} \\ \hline
NLP-IP    & 1 (Ref.)    & 1.01e-8     & 25.0    & 76\\
NLP-SQP    & 1.001     & 3.31e-8    & 14.8    & 62\\
LTV-NN    & 1.017    & 3.51e-3    & 0.71    & 94\\
LTI-NN    & 1.046    & 7.86e-3    & 0.67    & 42\\ 
LTV-Poly    & -    & -    & -    & 0\\ 
LTI-Poly    & -    & -    & -    & 0\\
LTV-NT    & 0.988    & $-$4.35e-3    & 0.70    & 100\\ \hline
\end{tabular}
\end{table}

\section{Conclusion}
\label{sec:conclusion}

This paper presented a learning-based linear time-varying (LTV) Koopman framework for constrained nonlinear optimal control. Motivated by the challenges of representing controlled nonlinear dynamics and nonlinear constraints within a computationally efficient optimal control framework, the proposed approach integrates Koopman-based modeling, neural-network learning, and convex optimization into a unified formulation.

Compared with conventional linear time-invariant Koopman formulations, the proposed LTV Koopman representation provides a more expressive description of controlled nonlinear dynamics while preserving a linear prediction structure suitable for efficient optimization. To support the use of the resulting model in optimal control, approximation errors and their propagation over extended prediction horizons were analyzed. Building upon these results, a neural-network-based learning framework was developed to jointly learn lifted dynamics and constraint representations from data, together with a stability-oriented regularization strategy to improve long-horizon prediction performance.

Using the learned Koopman model, the original constrained nonlinear optimal control problem was reformulated as a convex quadratic program. To account for modeling uncertainty, an error-aware constraint-tightening strategy was further introduced to improve the reliability of constraint enforcement. Numerical studies on a low-thrust orbit transfer problem and a six-degree-of-freedom Mars entry guidance problem demonstrated that the proposed framework provides improved prediction accuracy and constraint representation compared with conventional LTI Koopman formulations, while generating control solutions that remain competitive with nonlinear programming approaches at substantially lower computational cost.

Overall, the proposed framework establishes a unified connection between Koopman-based modeling, data-driven learning, and constrained optimal control. By combining an LTV Koopman representation, joint learning of dynamics and constraints, and error-aware convex optimization, it provides a practical and computationally efficient approach for solving constrained nonlinear optimal control problems.

}

\section*{Acknowledgments}
This project is supported in part by NSF 2334995 and NSF-2523830. An AI tool has been used exclusively for grammar and language editing.

\bibliography{reflib}

\end{document}